\documentclass[12pt,leqno]{article}
\usepackage{amsfonts}
\usepackage{amsmath, amsthm, amsfonts, amssymb, color}
\usepackage{mathrsfs}
\usepackage{color}
\theoremstyle{definition}

\newcommand{\scr}[1]{\mathscr #1}
\definecolor{wco}{rgb}{0.5,0.2,0.3}

\numberwithin{equation}{section} \theoremstyle{remark}

\newcommand{\ua}{\uparrow}

\title{{\bf  Space-Distribution  PDEs for Path Independent Additive Functionals of  McKean-Vlasov SDEs}\footnote{Supported in
 part by  NNSFC (11771326, 11431014).} }
\author{
{\bf Panpan Ren$^{b)}$,    Feng-Yu Wang$^{a,b)}$  }\\
\footnotesize{$^{a)}$ Center for Applied Mathematics, Tianjin University, Tianjin 300072, China}\\
 \footnotesize{$^{b)}$ Department of Mathematics,
Swansea University, Singleton Park, SA2 8PP, United Kingdom}\\
\footnotesize{  673788@swansea.ac.uk, wangfy@tju.edu.cn, F.-Y.Wang@swansea.ac.uk}}
\begin{document}
\allowdisplaybreaks
\def\R{\mathbb R}  \def\ff{\frac} \def\ss{\sqrt} \def\B{\mathbf
B}
\def\N{\mathbb N} \def\kk{\kappa} \def\m{{\bf m}}
\def\ee{\varepsilon}\def\ddd{D^*}
\def\dd{\delta} \def\DD{\Delta} \def\vv{\varepsilon} \def\rr{\rho}
\def\<{\langle} \def\>{\rangle} \def\GG{\Gamma} \def\gg{\gamma}
  \def\nn{\nabla} \def\pp{\partial} \def\E{\mathbb E}
\def\d{\text{\rm{d}}} \def\bb{\beta} \def\aa{\alpha} \def\D{\scr D}
  \def\si{\sigma} \def\ess{\text{\rm{ess}}}
\def\beg{\begin} \def\beq{\begin{equation}}  \def\F{\scr F}
\def\Ric{\mathcal Ric} \def\Hess{\text{\rm{Hess}}}
\def\e{\text{\rm{e}}} \def\ua{\underline a} \def\OO{\Omega}  \def\oo{\omega}
 \def\tt{\tilde}
\def\cut{\text{\rm{cut}}} \def\P{\mathbb P} \def\ifn{I_n(f^{\bigotimes n})}
\def\C{\scr C}      \def\aaa{\mathbf{r}}     \def\r{r}
\def\gap{\text{\rm{gap}}} \def\prr{\pi_{{\bf m},\varrho}}  \def\r{\mathbf r}
\def\Z{\mathbb Z} \def\vrr{\varrho} \def\ll{\lambda}
\def\L{\scr L}\def\Tt{\tt} \def\TT{\tt}\def\II{\mathbb I}
\def\i{{\rm in}}\def\Sect{{\rm Sect}}  \def\H{\mathbb H}
\def\M{\scr M}\def\Q{\mathbb Q} \def\texto{\text{o}} \def\LL{\Lambda}
\def\Rank{{\rm Rank}} \def\B{\scr B} \def\i{{\rm i}} \def\HR{\hat{\R}^d}
\def\to{\rightarrow}\def\l{\ell}\def\iint{\int}
\def\EE{\scr E}\def\Cut{{\rm Cut}}\def\W{\mathbb W}
\def\A{\scr A} \def\Lip{{\rm Lip}}\def\S{\mathbb S}
\def\BB{\scr B}\def\Ent{{\rm Ent}} \def\i{{\rm i}}\def\itparallel{{\it\parallel}}
\def\g{{\mathbf g}}\def\Sect{{\mathcal Sec}}\def\T{\mathcal T}
\def\f{\mathbf f} \def\g{\mathbf g}
\maketitle

\begin{abstract} Let $\scr P_2(\R^d)$ be the space of probability measures on $\R^d$ with finite second moment. The path independence of additive functionals of McKean-Vlasov SDEs is characterized by   PDEs on the product space $\R^d\times \scr P_2(\R^d)$ equipped with the usual derivative in space variable and Lions' derivative  in distribution. These PDEs are solved by using probabilistic arguments developed from \cite{LP}. As consequence, the path independence of   Girsanov transformations are identified with   nonlinear PDEs on  $\R^d\times \scr P_2(\R^d)$ whose solutions are given by probabilistic arguments as well. In particular,  the corresponding results on the Girsanov transformation killing the drift term derived earlier for the classical SDEs are recovered as special situations.
\end{abstract} \noindent
 \vskip 2cm

\section{Introduction}

In recent years, McKean-Vlasov stochastic differential equations (SDEs), also called distribution dependent or mean field SDEs, have received increasing attentions for their theoretically importance in characterizing  non-linear Fokker-Planck equations from physics. On the other hand, SDEs have been developed as crucial mathematical tools modelling economic and finance systems. In the real world, the evolution of these systems is not only driven by  micro actions (drift and noise), but also relies on  the macro environment (in mathematics, distribution of the systems). So,  it is reasonable to characterize economic and finance systems by using distribution dependent SDEs.

Let $\scr P(\R^d)$ be the space of all probability measures on $\R^d$, and let
 $$\scr P_2(\R^d)=\bigg\{\mu\in \scr P(\R^d):\ \mu(|\cdot|^2):=\int_{\R^d} |x|^2\mu(\d x)<\infty\bigg\}.$$ Then $\scr P_2(\R^d)$ is a Polish space under the Wasserstein distance
 $$\W_2(\mu,\nu):= \inf_{\pi\in \C(\mu,\nu)} \bigg(\int_{\R^d\times\R^d} |x-y|^2\pi(\d x,\d y)\bigg)^{\ff 1 2},\ \ \mu,\nu\in \scr P(\R^d),$$
 where $\C(\mu,\nu)$ is the set of couplings for $\mu$ and $\nu$; that is, $\pi\in\C(\mu,\nu)$ is a probability measure on $\R^d\times\R^d$ such that
 $\pi(\cdot\times\R^d)=\mu$ and $\pi(\R^d\times\cdot)=\nu$.

 Let $W_t$ be an $m$-dimensional Brownian motion on a standard filtered probability space $(\OO,\F,\{\F_t\}_{t\ge 0},\P)$, and denote by $\L_\xi$ the distribution of a random variable $\xi$ on $\R^d$.
Consider the following McKean-Vlasov SDE on $\R^d$:
\beq\label{E1} \d X_t=  b(t, X_t,\L_{X_t})\d t + \si(t,X_t,\L_{X_t}) \d W_t,\end{equation}
where
$$\si: [0,\infty)\times \R^d\times \scr P_2(\R^d)\to \R^{d\otimes m},\ \ b: [0,\infty)\times \R^d\times \scr P_2(\R^d)\to \R^d$$
are continuous such that for some increasing function $K: [0,\infty)\to [0,\infty)$ there holds
\beq\label{SPP1}
\beg{split} &|  b(t,x,\mu)- b(t,y,\nu)|+\|\si(t,x,\mu)-\si(t,y,\nu)\|_{HS}\\
&\le K(t)\big(|x-y|+\W_2(\mu,\nu)\big),\ \ t\ge 0, x,y\in \R^d, \mu,\nu\in \scr P_2(\R^d)
\end{split}\end{equation}
 and
 \beq\label{SPPn}
\|\si(t,{\bf 0},\dd_{{\bf 0}})\|_{HS}+ |b(t,{\bf 0},\dd_{{\bf 0}})|\le K(t),\ \ t\ge 0,
\end{equation}
where $\dd_{{\bf 0}}$ is the Dirac measure at ${\bf 0}\in\R^d$. For any $t\ge 0$, let $L^2(\OO\to\R^d,\F_t,\P)$ be the class of $\F_t$-measurable square integrable random variables on $\R^d$. By \eqref{SPP1} and \eqref{SPPn},
for any $s\ge 0$ and $X_s\in L^2(\OO\to\R^d,\F_s,\P)$, \eqref{E1} has a unique solution $(X_t)_{t\ge s}$ with

\beq\label{MM}
\sup_{t\in [s,T]} \E|X_t|^2<\infty,\ \ T\ge s.
\end{equation}
See  \cite{W18} for more results on gradient estimates and Harnack inequalities of the associated nonlinear semigroup, and \cite{18HW,MV} and references within for   the existence and uniqueness under weaker conditions.

In this paper, we aim to characterize the path independence of the additive functional
\beq\label{YPP2}
A_{s,t}^{\f,\g}:= \int_s^t \f(r,X_r,\L_{X_r}) \d r + \int_s^t\<\g(r,X_r,\L_{X_r}),\d W_r\>,\ \ 0\le s\le t,
\end{equation}
where $$\f: [0,\infty)\times \R^d\times \scr P_2(\R^d)\to \R,\ \ \g: [0,\infty)\times \R^d\times \scr P_2(\R^d)\to \R^m$$
are continuous, so that $A_{s,t}^{\f,\g}$ for $t\ge s$ is a well-defined local semi-martingale.

\beg{defn} The additive functional  $(A_{s,t}^{\f,\g})_{t\ge s}$ is called path independent, if there exists a measurable function
$$V: [0,\infty)\times \R^d\times\scr P_2(\R^d)\to\R$$ such that
\beq\label{YPP0} A_{s,t}^{\f,\g} = V(t,X_t, \L_{X_t})- V(s, X_s, \L_{X_s}),\ \ 0\le s\le t.\end{equation} \end{defn}

The motivation of the study comes from mathematical statement of equilibrium financial market. In their seminal paper \cite{BS} Black and Scholes described the price dynamics (or the wealth growth) by using SDEs under a so-called real world probability  measure. But for an equilibrium financial market there exists a so-called risk neutral   measure having a path independent density with respect to the real world probability, see \cite{HC}. That is,   under the risk neutral  measure the solution of \eqref{E1} becomes a martingale, and  the density of the neutral probability with respect to the real world one depends only on the initial and current states but not those in between.

For instance, let $\f=\ff 1 2|\g|^2$. Then $A_{s,t}^{\f,\g}$ becomes
\beq\label{APPg} A_{s,t}^\g:= \ff 1 2\int_s^t |\g(r,X_r,\L_{X_r})|^2\d r + \int_s^t\<\g(r,X_r,\L_{X_r}),\d W_r\>,\ \ 0\le s\le t. \end{equation}
By the Girsanov theorem, when \beq\label{Agg2}\E\e^{\ff 1 2\int_s^t |\g(r,X_r,\L_{X_r})|^2\d r}<\infty,\end{equation}  $\d \Q_{s,t}^\g:= \e^{-A_{s,t}^\g}\d\P$  is a probability measure.  So, to adopt $\Q_{s,t}^\g$ as a risk neutral measure, we need to verify   the path independence of the additive functional $A_{s,t}^\g$ in the sense of \eqref{YPP0}. In particular, when
\beq\label{BB} b=\si \tt b\ \text{for\ some\ measurable\ } \tt b:  [0,\infty)\times \R^d\times \scr P_2(\R^d)\to \R^m,\end{equation}
and  \eqref{Agg2} holds for $\g:=\tt b$,   let
\beq\label{YPP1}
A_{s,t}:= \ff 1 2 \int_s^t |\tt b  (r,X_r,\L_{X_r})|^2 \d r + \int_s^t\<\tt b(r,X_r,\L_{X_r}),\d W_r\>,\ \ 0\le s\le t.
\end{equation}
Then  $\d \Q_{s,t}:= \e^{-A_{s,t}}\d\P$  is a probability measure such that
$$\tt W_r:= W_r+ \int_s^r  \tt b(u,X_u,\L_{X_u})\d u,\ \ r\in [s,t]$$ is an $m$-dimensional Brownian motion, and hence
$$X_r= X_s+\int_s^r \si(u,X_u,\L_{X_u})\d \tt W_u,\ \ r\in [s,t]$$ is  a $\Q_{s,t}$-martingale as required for an equilibrium financial market.   We   would like to investigate the path independence of the additive functional $A_{s,t}$ such that $\Q_{s,t}$ is a risk neutral measure.

In general, to characterize the path independence of $A_{s,t}^{\f,\g}$ using  PDEs on $\R^d\times\scr P_2(\R^d)$, we shall need that   $\L_{X_t} (t>s)$ has a full support on $\R^d$. To this end, we will assume  the   H\"ormander condition.
To state this condition, for any $\mu_\cdot\in C([0,\infty);\scr P_2(\R^d))$, let
 \beg{align*} & U_j^\mu(t,x)= \sum_{i=1}^d \si_{ij}(t,x,\mu_t) \pp_{x_i},\ \  1\le j\le m,\\
& U_0^\mu(t,x)= \sum_{i=1}^d \Big[b_i(t,x,\mu_t) - \ff 1 2 \sum_{j=1}^m\sum_{k=1}^d \big\{\si_{kj}\pp_{x_k}\si_{ij}\big\}(t,x,\mu_t)\Big]\pp_{x_i}.\end{align*}
Then  for  $\mu_t:=\L_{X_t}$, the SDE \eqref{E1} reduces to
\beg{align*} \d X_t &=b(t,X_t,\mu_t)\d t + \si(t,X_t,\mu_t)\d W_t\\
&= U_0^\mu(t,X_t)\d t +\sum_{j=1}^m U_j^\mu(t,X_t)\circ \d W^j_t,\end{align*}
where $\circ \d W^j_t$ is the Stratonovich differential with respect to the $j$-th component of $W_t$. Let
 $\scr U_0^\mu=\{U_j^\mu: 1\le j\le m\}$ and $$\scr U_{l}^\mu=\scr U_0^\mu\cup\big\{{\rm Lie\ brackets\ up\ to\ order}\ l {\rm\ for\ vector\ fields}\  U_j^\mu: 0\le j\le m \big\},\ \ l\in \mathbb N.$$
The  H\"ormander condition \cite{HOR} is stated as follows.

\beg{enumerate}\item[$(H)$] For any $\mu_\cdot\in C([0,\infty);\scr P_2(\R^d))$, there exists $l\in \mathbb Z_+$ such that the vector fields $\{U_j: 0\le j\le m\}$ are $C^l$-smooth and the    family of vector fields $\scr U_l^\mu$ span    the tangent space. \end{enumerate}
By the Harnack inequality for hypoelliptic parabolic equations,  the  condition $(H)$  implies that for any $s\ge 0$ and
$X_s\in L^2(\OO\to\R^d,\F_s,\P)$, the the distribution  $\L_{X_t} $ for $t>s$ has  has full support on $\R^d$, see the proof of Lemma \ref{LNN} below for details.

In Section 2, we will characterize  the path independence of $A_{s,t}^{\f,\g}$ using PDEs on $\R^d\times\scr P_2(\R^d)$, see Theorem \ref{T1.1} below for details. Following the idea of \cite{LP}, such type PDEs are solved using solutions of an associated SDE, see Theorem \ref{T1.2} for details. As a consequence,   the path independence of $A_{s.t}^\g$ in \eqref{APPg}  and $A_{s,t}$ in \eqref{YPP1} is identified with nonlinear PDEs on $\R^d\times \scr P_2(\R^d)$, see Corollaries \ref{T1.4} and  \ref{C1.3} below. When  the SDE is distribution independent, i.e. $b(t,x,\mu)$ and $\si(t,x,\mu)$ do not depend on $\mu\in \scr P_2(\R^d)$, Corollary \ref{C1.3} recovers the corresponding existing results derived in \cite{TWW, WW, YW},  see also \cite{QW,WW1} for extensions to SDEs with jumps and semi-linear SPDEs.
  Finally, complete proofs of these results are presented in Section 3.

\section{Main results  }

To state our results, we first recall  the definition of $L$-derivative for functions on $\scr P_2(\R^d)$, which was introduced by P.-L. Lions in his  lectures \cite{Card} at College de France, see also \cite{LP, HSS}. In the following we introduce a straightforward definition without  using abstract probability spaces as  in previous references.   Let $\pp_t$ denote the partial differential in time parameter $t\ge 0$, $\pp_x$ or $\pp_y$ the gradient operator in variables $x$ or $y\in \R^d$, and $\pp^2_x$ the Hessian operator in $x\in \R^d$.
Let ${\rm Id}:\mathbb{R}^{d}\rightarrow\mathbb{R}^{d}$ be the identity map, i.e.  ${\rm Id}(x)=x$ for $x\in \R^d$. It is easy to see that for any $\mu\in \scr P_2(\R^d)$ and $\phi\in L^2(\R^d\to\R^d,\mu)$, we have $\mu\circ ({\rm Id}+\phi)^{-1} \in \scr P_2(\R^d).$

\beg{defn}\label{defn1}  Let   $T\in (0,\infty]$,  and   set $[0,T]= [0,\infty)$ when $T=\infty$.
\beg{enumerate} \item[$(1)$] A function $f: \scr P_2(\R^d)\to \R$ is called  $L$-differentiable at $\mu\in \scr P_2(\R^d)$, if  the functional
$$L^2(\R^d\to \R^d,\mu)\ni \phi \mapsto f(\mu\circ({\rm Id} +\phi)^{-1})$$ is Fr\'echet differentiable at ${\bf0}\in L^2(\R^d\to \R^d,\mu)$;
that is, there exists $($hence, unique$)$ $\xi\in L^2(\R^d\to\R^d,\mu)$ such that
\beq\label{*D1}\lim_{\mu(|\phi|^2)\to 0} \ff{f(\mu\circ ({\rm Id}+\phi)^{-1})-f(\mu)- \mu(\<\xi,\phi\>)}{\ss{\mu(|\phi|^2)}}= 0.\end{equation}
In this case, we denote $\pp_\mu f(\mu)=\xi$ and call it the $L$-derivative of $f$ at $\mu$.

 \item[$(2)$]  A function $f: \scr P_2(\R^d)\to \R$ is called $L$-differentiable on $\scr P_2(\R^d)$ if  the $L$-derivative $\pp_\mu f(\mu)$ exists for all $\mu\in\scr P_2(\R^d)$. If moreover
 $(\pp_\mu f(\mu))(y)$ has a version  differentiable in $y\in \R^d$ such that $(\pp_\mu f(\mu))(y)$ and $\pp_y (\pp_\mu f(\mu))(y)$ are jointly continuous in $(\mu,y)\in \scr P_2(\R^d)\times\R^d$, we denote $f\in C^{(1,1)}(\scr P_2(\R^d))$.

 \item[$(3)$] A function $f: [0,T]\times\R^d\times\scr P_2(\R^d)\to \R$ is said to be in the class $C^{1,2,(1,1)}([0,T]\times\R^d\times\scr P_2(\R^d))$, if the derivatives
$$\pp_t f(t,x,\mu), \ \pp_x f(t,x,\mu), \ \pp_x^2 f(t,x,\mu), \ \pp_\mu f(t,x,\mu)(y),\  \pp_y\pp_\mu f(t,x,\mu)(y)$$ exist
 and are jointly continuous in the corresponding arguments $(t,x,\mu)$ or $(t,x,\mu,y)$. If $f\in C^{1,2,(1,1)}([0,T]\times\R^d\times\scr P_2(\R^d))$ with all these  derivatives   bounded on $[0,T]\times\R^d\times\scr P_2(\R^d)$, we denote $f\in C^{1,2,(1,1)}_b([0,T]\times\R^d\times\scr P_2(\R^d))$.

  \item[$(4)$]  Finally, we write  $f\in \scr C([0,\infty)\times\R^d\times\scr P_2(\R^d))$, if $f\in C^{1,2,(1,1)}([0,T]\times\R^d\times\scr P_2(\R^d))$ and  the function $$ (t,x,\mu)\mapsto  \int_{\R^d} \big\{\|\pp_y\pp_\mu f\|+\|\pp_\mu f\|^2\big\}(t,x,\mu)(y)\mu(\d y)$$
  is locally bounded, i.e. it is bounded on compact subsets of   $[0,T]\times\R^d\times\scr P_2(\R^d).$
\end{enumerate}
\end{defn}

For readers' understanding of the $L$-derivative, we present below an example for  a  class of functions inducing the Borel $\si$-algebra on $\scr P_2(\R^d)$. See \cite[Example 2.2]{LP} for  concrete choices of $F$ and $h_i$.

\beg{exa} Let $n\in \mathbb N, \{h_i\}_{1\le i\le n}\subset C^2(\R^d)$ with  $\|\pp_x^2h_i\|_\infty<\infty$ and let $F\in C^1(\R^n)$. Then
the function $$\scr P_2(\R^d)\ni\mu\mapsto f(\mu):= F(\mu(h_1),\cdots, \mu(h_n)) $$ is in $C^{(1,1)}(\scr P_2(\R^d))$ with
$$\pp_\mu f(\mu)(y)= \sum_{i=1}^n (\pp_i F)(\mu(h_1),\cdots, \mu(h_n)) \pp_yh_i(y),\ \ \mu\in \scr P_2(\R^d), y\in \R^d.$$\end{exa}
\beg{proof}
By the chain rule it suffices to prove for $f(\mu):= \mu(h_1)$, i.e. $n=1$ and $F(r)=r$.  Since $\|\pp_x^2h_i\|_\infty<\infty$, there exists a constant $C>0$ such that
$$  |h_i(x)|+|\pp_x h_i(x)|^2 \le C(1+|x|^2),\ \ x\in \R^d,$$
so that $h_1\in L^1(\mu)$ and $\pp_xh_1\in L^2(\R^d\rightarrow\R^d, \mu)$ for $\mu\in \scr P_2(\R^d)$.
 Then, for any $\phi\in L^2(\R^d\rightarrow\R^d,\mu)$, by Taylor's expansion we have
\beg{align*} &\lim_{\|\phi\|_{L^2(\mu)}\to 0} \ff{|f(\mu\circ({\rm Id}+\phi))-f(\mu)-\mu(\<\pp h_1, \phi\>)|}{\|\phi\|_{L^2(\mu)}}\\
&= \lim_{\|\phi\|_{L^2(\mu)}\to 0} \ff 1 {\|\phi\|_{L^2(\mu)}}  \bigg| \int_{\R^d}\big\{h_1(x+\phi(x))-h_1(x)-\<\pp_x h_1(x),\phi(x)\>\big\}\mu(\d x)\bigg|\\
&\le   \lim_{\|\phi\|_{L^2(\mu)}\to 0} \ff{\|\pp_x^2h_1\|_\infty}{ 2\|\phi\|_{L^2(\mu)} } \int_{\R^d}|\phi(x)|^2\mu(\d x)  \le \lim_{\|\phi\|_{L^2(\mu)}\to 0}  \|\pp_x^2h_1\|_\infty \|\phi\|_{L^2(\mu)}= 0.\end{align*}So,  by definition, $\pp_\mu f(\mu)(y)= \pp_y h_1(y).$\end{proof}

Let us explain that  the above definition of $L$-derivative coincides with the Wasserstein derivative introduced by P.-L. Lions using probability spaces.  Given $\mu \in \scr P_2(\R^d)$, let $(\tt \OO,\tt \F,\tt \P)=(\R^d,\scr B(\R^d),\mu)$ and  $X={\rm Id}$. Then $X$ is a random variable with $\L_X|_{\tt \P}=\mu$.  For any square integrable random variable $Y$,   we have $\phi:=Y\in L^2(\R^d\to\R^d,\mu)$. Moreover, since $X={\rm Id}$,   for any $A\in\scr B(\R^d)$,
\begin{equation*}
\begin{split}
 (\mu\circ ({\rm Id}+\phi)^{-1})(A)&=\mu(\{x:({\rm Id}+\phi)(x)\in A\})
 =\mu(\{x:x+\phi(x)\in A\})\\
 &=
 \tt P(\{x:X(x)+Y(x)\in A\})=\tt P(X+Y\in A)\\
 &=\L_{X+Y}|_{\tt \P}(A).
\end{split}
\end{equation*}
So,  $(\mu\circ ({\rm Id} +\phi)^{-1})=(\L_{X+Y}|_{\tt \P})$, and \eqref{*D1} means that
$$L^2(\tt\OO\to\R^d,\tt\P)\ni Y\mapsto f(\L_{X+Y}|_{\tt\P})$$ is Fr\'echet differentiable with derivative $\pp_\mu f(\mu):=\xi$, which coincides with \cite[Definition 6.1]{Card} given by P.-L. Lions. Note that the atomless restriction on the probability space therein is to ensure the existence of a random variable with distribution   $\mu$.
It is crucial that (see \cite[Proposition A.2]{HSS}) the definition of $\pp_\mu f(\mu)\in L^2(\mathbb{R}^{d}\rightarrow\mathbb{R}^{d}, \mu)$ dose not depend on the choice of probability space
$(\tt \OO,\tt \F,\tt \P)$ and random variable $X$ with $\L_X|_{\tt \P}=\mu$. So, in particular, we may take the above specific choice $(\tt \OO,\tt \F,\tt \P)=(\R^d,\scr B(\R^d),\mu)$ and $X={\rm Id}$.

 The following differential operator on $[0,\infty) \times \R^d\times\scr P_2(\R^d)$ associated with the SDE \eqref{E1} has been introduced in \cite{LP}: for any $V\in C^{1,2,(1,1)}([0,\infty)\times\R^d\times\scr P_2(\R^d))$ and $(t,x,\mu)\in [0,\infty)\times  \R^d\times \scr P_2(\R^d)$, let
\beg{equation} \label{XTT}\begin{split}{\bf L}_{\si,b} V(t,x,\mu)&=  \ff 1 2 {\rm tr}\big(\si\si^*  \pp_{x}^2 V)(t,x,\mu) +\<b, \pp_x V\>  (t,x,\mu)\\
&\quad+\int_{\R^d} \Big[\ff 1 2 {\rm tr}\big\{(\si\si^*)(t,y,\mu)  \pp_{y}\pp_\mu V (t,x,\mu)(y)\big\} \\
&\quad+\big\<b(t,y,\mu) ,  \pp_{\mu} V(t,x,\mu)(y)\big\>\Big]\mu(\d y).\end{split}\end{equation}
Our first result  is the following characterization on the path independence of the functional $A_{s,t}^{\f,\g}$ in \eqref{YPP2}.

\beg{thm} \label{T1.1} Assume that $\sigma$ and $b$ satisfy \eqref{SPP1} and \eqref{SPPn} for some locally bounded function $K$.  Let $T>0$, $\f\in C([0,T]\times \R^d\times\scr P_2(\R^d))$
and $\g\in C([0,T]\times \R^d\times\scr P_2(\R^d)\to \R^m)$.
\beg{enumerate} \item[$(1)$] If $V\in \scr C ([0,T]\times\R^d\times\scr P_2(\R^d))$ solves the PDE
 \beq\label{ATT0} \beg{cases}(\pp_t+ {\bf L}_{\si,b}) V (t,x,\mu)=\f(t,x,\mu),\\
 (\si^*\pp_x V)(t,x,\mu)=\g(t,x,\mu),  \end{cases} \ \ t\in [0, T], x\in \R^d, \mu\in \scr P_2(\R^d),
\end{equation} then the path independent property $\eqref{YPP0}$ holds.
\item[$(2)$] Conversely, if $(H)$ holds, then  the path independence of $A_{s,t}^{\f,\g}$  in the sense of  $\eqref{YPP0}$  for some $V\in \scr C ([0,T]\times\R^d\times\scr P_2(\R^d))$ implies  $\eqref{ATT0}.$\end{enumerate}
\end{thm}

To provide a class of $(\f,\g)$ such that the additive functional $A_{s,t}^{\f,\g}$ is path independent in the sense of \eqref{YPP0}, we adopt the idea of \cite{LP} to solve the PDE \eqref{ATT0} using  an SDE accompying with  \eqref{E1}. To state this accompying SDE, for any $\mu\in \scr P_2(\R^d)$ and $s\ge 0$, let $(X_{s,t}^\mu)_{t\ge s}$ solve \eqref{E1} from time $s$ with $\scr L_{X_{s,s}^\mu}=\mu$. Let
\beq\label{SM0}
 P_{s,t}^*\mu = \L_{X_{s,t}^\mu}, \ \ \ t\ge s,\ \ \mu\in \scr P_2(\R^d).
\end{equation}
As shown in \cite{W18} that  $P_{s,t}^*$ is a nonlinear semigroup satisfying
\beq\label{SM}
P_{t,r}^*P_{s,t}^*= P_{s,r}^*,\ \ \ 0\le s\le t\le r.
\end{equation} Now, for any $x\in \R^d, \mu\in \scr P_2(\R^d)$ and $s\ge 0$, let $(X_{s,t}^{x,\mu})_{t\ge s}$ solve the SDE
\beq\label{E1'*} \d X_{s,t}^{x,\mu}= b(t, X_{s,t}^{x,\mu}, P_{s,t}^*\mu)\d t +\si(t, X_{s,t}^{x,\mu}, P_{s,t}^*\mu)\d W_t,\ \ X_{s,s}^{x,\mu}=x.
\end{equation}
We have the following result.

\beg{thm}\label{T1.2} Assume that $b_i,\si_{ij},\f\in C_b^{1,2,(1,1)}([0,T] \times \R^d\times\scr P_2(\R^d)),\ 1\le i\le d, 1\le j\le m$.   Then $V\in C_b^{1,2,(1,1)}([0,T]\times \R^d\times \scr P_2(\R^d))$ solves the first PDE in  $\eqref{ATT0}$  if  there exists $\Phi\in C_b^{2,(1,1)}(\R^d\times \scr P_2(\R^d))$ such that
$$ V(t,x,\mu)= \E \bigg(\Phi(X_{t,T}^{x,\mu},P_{t,T}^*\mu)  -\int_t^T \f(r,X_{t,r}^{x,\mu}, P_{t,r}^*\mu)\d r\bigg),  \ \ t\in [0, T], x\in \R^d, \mu\in \scr P_2(\R^d).
 $$ Consequently, for given $V\in C_b^{1,2,(1,1)}([0,T]\times\R^d\times\scr P_2(\R^d))$,  $A_{s,t}^{\f,\g}$ is path independent  in the sense of $\eqref{YPP0}$  if
 $$\g(t,x,\mu)= \si^* \pp_x \E \bigg(V(T,X_{t,T}^{x,\mu},P_{t,T}^*\mu)  -\int_t^T \f(r,X_{t,r}^{x,\mu}, P_{t,r}^*\mu)\d r\bigg)$$ for all $(t,x,\mu)\in[0, T]\times \R^d\times\scr P_2(\R^d),$
  and the inverse holds under assumption $(H)$.\end{thm}

Next,  we consider  $f:= \ff 1 {2\bb} |\g|^2$  for a constant $\bb\ne 0$. Then  the additive functional $A_{s,t}^{\f,\g}$ reduces to
\beq\label{APPg3} A_{s,t}^{\g;\bb}:= \ff 1 {2\bb} \int_s^t |\g(r,X_r,\L_{X_r})|^2 \d r + \int_s^t\<\g(r,X_r,\L_{X_r}),\d W_r\>,\ \ 0\le s\le t.
\end{equation} This covers $A_{s,t}^{\g}$ in \eqref{APPg} for $\bb=1.$
As a consequence of Theorems \ref{T1.1} and \ref{T1.2}, we have the following result on  the path independence of $A_{s,t}^{\g;\bb}$ and the corresponding nonlinear PDE:
\beq\label{NPDE} (\pp_t +{\bf L}_{\si,b})V(t,x,\mu)= \ff 1{2\bb} |\si^*\pp_x V|^2(t,x,\mu),\ \ (t,x,\mu)\in [0,T]\times\R^d\times\scr P_2(\R^d).\end{equation}

\beg{cor}\label{T1.4} Assume that $\sigma$ and $b$ satisfy \eqref{SPP1} and \eqref{SPPn} for some locally bounded function $K$.  Let $T>0$ and  $0\ne \bb\in\R$.
\beg{enumerate} \item[$(1)$] If $V\in \scr C ([0,T]\times\R^d\times\scr P_2(\R^d))$ solves the nonlinear PDE $\eqref{NPDE}$, then $A_{s,t}^{\g;\bb}$ with  $\g:= \si^*\pp_x V$ is path independent in the sense of  $\eqref{YPP0}.$
Conversely, under assumption $(H)$, for any $\g\in C([0,T]\times \R^d\times\scr P_2(\R^d)\to \R^d)$, the path independence of $A_{s,t}^{\g;\bb}$ in the sense of  $\eqref{YPP0}$  for some $V\in \scr C ([0,T]\times\R^d\times\scr P_2(\R^d))$ implies that $\g=\si^*\pp_xV$ and $V$ solves $\eqref{NPDE}$.
 \item[$(2)$] Let $b_i,\si_{ij} \in C_b^{1,2,(1,1)}([0,T] \times \R^d\times\scr P_2(\R^d))$ for $1\le i\le d, 1\le j\le m$. For any
$\Phi\in C^{2,(1,1)}_b(\R^d\times\scr P_2(\R^d))$ with $\inf\Phi>0$,
\beq\label{TTY0} V(t,x,\mu):= -\bb \log \big\{\E \Phi(\tt X_{t,T}^{x,\mu}, \tt \mu_{t,T})\big\},\ \ (t,x,\mu)\in [0,T]\times\R^d\times\scr P_2(\R^d)\end{equation} is the unique solution to the nonlinear PDE
$\eqref{NPDE}$ with
$$V(T,x,\mu)= -\bb \log \Phi(x,\mu),\ \ \ (x,\mu)\in \R^d\times \scr P_2(\R^d).$$\end{enumerate}
\end{cor}

Finally, we consider the path independence of   the functional $A_{s,t}$ in \eqref{YPP1}.
Let
\beq\label{SIG}\beg{split}  {\bf L}_{\si} V(t,x,\mu)=  &\ff 1 2 {\rm tr}\big(\si\si^*  \pp_{x}^2 V)(t,x,\mu) +\ff 1 2|\si^* \pp_x V|^2  (t,x,\mu)\\
&+ \int_{\R^d} \Big[\ff 1 2{\rm tr}\big\{(\si\si^*)(t,y,\mu)  \pp_{y}\pp_\mu V (t,x,\mu)(y)\big\}\\
&\qquad  +\big\<(\si\si^*\pp_yV )(t,y,\mu) ,  \pp_{\mu} V(t,x,\mu)(y)\big\>\Big]\mu(\d y).\end{split}\end{equation}

\beg{cor} \label{C1.3} Assume that $\sigma$ and $b$ satisfy \eqref{SPP1} and \eqref{SPPn} for some locally bounded function $K$. Let $T>0$.
\beg{enumerate} \item[$(1)$]  If $V\in \scr C ([0,T]\times\R^d\times\scr P_2(\R^d))$ solves the nonlinear PDE
 \beq\label{ATT2} \beg{cases}(\pp_t+ {\bf L}_{\si}) V (t,x,\mu)=0,\\
 b(t,x,\mu)= (\si\si^*\pp_x V)(t,x,\mu),  \end{cases}\ \ t\in [0, T], x\in \R^d, \mu\in \scr P_2(\R^d),  \end{equation}
then $\eqref{BB}$ holds for   $\tt b:=\si^*\pp_x V$ and  $A_{s,t}$ in $\eqref{YPP1} $ is path independent  in the sense of  $\eqref{YPP0}$. Conversely, under assumption $(H)$, if $\eqref{BB}$ holds for some $\tt b=\si^*\pp_x V$ such that  $A_{s,t}$ in $\eqref{YPP1} $ is path independent  in the sense of  $\eqref{YPP0}$ for some $V\in \scr C ([0,T]\times\R^d\times\scr P_2(\R^d))$, then $\eqref{ATT2}$ holds.
  \item[$(2)$] A function $V\in C_b^{1,2,(1,1)}([0,T]\times \R^d\times \scr P_2(\R^d))$ solves $\eqref{ATT2}$ if and only if there exists $\Phi\in C_b^{2,(1,1)}(\R^d\times \scr P_2(\R^d))$ with $\inf\Phi>0$ such that
$$\beg{cases}  V(t,x,\mu)= -\ff 1 2 \E \big\{\log \Phi(X_{t,T}^{x,\mu}, P_{t,T}^*\mu)\big\},\\
b(t,x,\mu)= (\si\si^*\pp_x V)(t,x,\mu),\  \ \ t\in [0, T], x\in \R^d, \mu\in \scr P_2(\R^d).\end{cases}
 $$\end{enumerate}
\end{cor}

Since $b(t,x,\mu)= (\si\si^*\pp_x V)(t,x,\mu)$ implies that both $X_{t,T}^{x,\mu}$ and $P_{t,T}^*\mu$ may depend on $V$, unlike Theorem \ref{T1.2} and Corollary \ref{T1.4}(2) providing solutions of \eqref{ATT0} and \eqref{NPDE} respectively, Corollary \ref{C1.3}(2)  only gives an alternative version of \eqref{ATT2} but not solutions.
To construct  a nontrivial solution  of \eqref{ATT2}, the nonlinear term
$$\int_{\R^d}  \big\<(\si\si^*\pp_yV )(t,y,\mu) ,  \pp_{\mu} V(t,x,\mu)(y)\big\> \mu(\d y) $$ in ${\bf L}_\si$ causes an essential difficulty. To overcome this difficulty, many other things   have to be treated. So, we would like to leave this problem to   a forthcoming paper.

 \section{Proofs }

We need the following It\^o's formula for distribution dependent functionals, see    \cite[Proposition 6.1]{LP} or \cite[Proposition A.8]{HSS} under stronger conditions on $\si$ and $f$.

\beg{lem}[It\^o's formula for distribution dependent
functional]\label{L1} For any $f\in \scr C ([0,\infty)\times\R^d\times\scr P_2(\R^d))$, $f(t,X_t,\L_{X_t})$ is a semi-martingale with
\begin{equation}\label{rpp}\d f(t,X_t, \L_{X_t}) = (\pp_t + {\bf L}_{\si,b}
)f(t,X_t,\L_{X_t}) \d t +\<(\si^*\pp_x f)(t,X_t, \L_{X_t}), \d
W_t\>,\end{equation}where ${\bf L}_{\si,b}$ is in \eqref{XTT}.
\end{lem}

\beg{proof} Let $\mu_t=\L_{X_t}$ and
$$\bar b(t,x)= b(t,x,\mu_t), \ \bar\si (t,x)=\si (t,x, \mu_t),\ \ \bar f(t,x)= f(t, x,\mu_t),\ \ t\ge 0, x\in \R^d.$$ Then $(X_t)_{t\ge 0}$ solves the classical SDE
$$\d X_t= \bar b (t,X_t)\d t+ \bar \si(t, X_t)\d W_t.$$
By the definition \ref{defn1} (4),  $f\in \scr C ([0,\infty)\times\R^d\times\scr P_2(\R^d))$ implies that $\bar f(t,x)$ is $C^2$-smooth in $x\in \R^d$. So,  if $\bar f(t,x)$ is $C^1$ in $t\ge 0$, we will be able to apply the classical It\^o's formula to derive
 \beg{equation*}
 \begin{split}\d f(t,X_t, \L_{X_t})&= \d \bar f(t,X_t) = \<(\bar\si^*\pp_x \bar f)(t,X_t), \d W_t\>\\
&\quad + \bigg\{ \pp_t\bar f  + \ff 1 2 \sum_{i,j=1}^d (\bar \si \bar \si^*)_{ij}  \pp_{x_i}\pp_{x_j} \bar f + \sum_{i=1}^d \bar b_i \pp_{x_i} \bar f\bigg\}(t,X_t)\d t\\
&= \<(\si^*\pp_x     f)(t,X_t,\L_{X_t}),   \d W_t\>+   (\pp_t\bar   f)(t,X_t)\d t\\
&\quad   + \bigg\{\ff 1 2 \sum_{i,j=1}^d (\si\si^*)_{ij}  \pp_{x_i}\pp_{x_j}   f + \sum_{i=1}^d b_i \pp_{x_i} f\bigg\}(t,X_t, \L_{X_t})\d t.
\end{split}
\end{equation*}
Therefore, to finish the proof, it suffices  to show that $\bar f(t,x)$ is differentiable in $t\ge 0$ and

\beq\label{QPP} \beg{split} (\pp_t \bar f)(t,x)& =\pp_t f(t,x,\nu)|_{\nu=\mu_t}+ \pp_t f(s,x,\mu_t)|_{s=t}\\
&=:(\pp_tf)(t,x,\mu_t)+(\pp_t^\mu f)(t,x,\mu_t),\end{split}
\end{equation}
where

\begin{equation}\label{QPP'}\beg{split}
(\pp_t^\mu f)(t,x,\mu_t):&=\ff 1 2\sum_{i,j=1}^d \int_{\R^d}  \Big[ (\si\si^*)_{ij}(t,y,\mu_t) \pp_{y_j} \{(\pp_\mu f)_i(t,x,\mu_t)\}(y)\Big] \mu_t(\d y)\\
&\quad+\sum_{i=1}^d \int_{\R^d} \Big[b_i(t,y,\mu_t) (\pp_\mu f)_i (t,x,\mu_t)\Big]\mu_t(\d y),
\end{split}
\end{equation}
is continuous in $(t,x)\in [0,\infty)\times \R^d$, since $f\in \scr C( [0,\infty)\times\R^d\times \scr P_2(\R^d))$ and
for any $T\in (0,\infty)$, $\{\mu_t:\ t\in [0,T]\}$ is a compact set in $\scr P_2(\R^d)$.  Below we prove \eqref{QPP} by two steps.

(a) According to \cite[Proposition A.6]{HSS}, if
\beq\label{QPP''} \E \int_0^T\big\{ |b(t,X_t,\mu_t)|^2+\| \si(t,X_t,\mu_t)\|_{HS}^4\big\}\d t<\infty,\end{equation} then for any $f\in \scr C( [0,\infty)\times\R^d\times \scr P_2(\R^d)),$
\beq\label{QPP5} \beg{split}  \bar f(t,x,\mu_{t+s})- f(t,x,\mu_t)
 = &\int_t^{t+s} \d r  \int_{\R^d} \Big[\ff 1 2  \sum_{i,j=1}^d (\si\si^*)_{ij} (r,x,\mu_r)\pp_{y_j}\{(\pp_\mu f)_i(r,y,\mu_r)\}(y)\\
 &\  +\sum_{i=1}^d b_i(r,y,\mu_r)(\pp_\mu f)_i(r,x,\mu_r)(y)\Big]\mu_r(\d y), \ \ s>0.\end{split}
\end{equation}
By conditions on  $b,\si$ and $f$, this implies \eqref{QPP}.

(b) In general, let $T>0$ be fixed. By \eqref{SPP1} and \eqref{SPPn} we have
\beq\label{MM2}|b(t,x,\mu_t)|^2+\|\si(t,x,\mu_t)\|_{HS}^2\le C \Big (1+|x|^2+\W_2^2(\mu_t,\delta_{\bf 0})\Big)     ,\ \ x\in\R^d, t\in [0,T]
\end{equation}
for some constant $C>0$. This, together with \eqref{MM}, implies
$$\E\int_0^T \big\{|b(t,X_t,\mu_t)|^2+ \|\si(t,X_t,\mu_t)\|_{HS}^2\big\}\d t<\infty.$$
So, to verify \eqref{QPP''}, we need to make approximations on $\si$. For any $k\in \mathbb N$, let
$$\phi_k(x)= \big(\{x_i\land k\}\lor\{-k\}\big)_{1\le i\le d},\ \ x\in \R^d.$$
Let $\si^{(k)}(t,x,\mu)= \si(t, \phi_k(x),\mu),$ and let $X_t^{(k)} $ solve the SDE
\beq\label{E11}   \d X_t^{(k)}= b(t, X_t^{(k)},\L_{X_t^{(k)}})\d t + \si^{(k)}(t,X_t^{(k)},\L_{X_t^{(k)}}) \d W_t,\ \ X_0^{(k)}=X_0.
\end{equation}
Then as explained in (a), $\mu_t^{(k)}:= \L_{X_t^{(k)}}$ satisfies
\beq\label{QPP5'} \beg{split}   & f(t,x, \mu_{t+s}^{(k)})- f(t,x,\mu_t^{(k)})\\
& =   \int_t^{t+s} \d r \int_{\R^d} \Big[\ff 1 2  \sum_{i,j=1}^d (\si^{(k)}(\si^{(k)})^*)_{ij} (r,x,\mu_r^{(k)})\pp_{y_j}\{(\pp_\mu f)_i(r,y,\mu_r^{(k)})\}(y)\\
 & \qquad +\sum_{i=1}^d b_i(r,y,\mu_r^{(k)})(\pp_\mu f)_i(r,x,\mu_r^{(k)})(y)\Big]\mu_r^{(k)}(\d y),\ \ s>0.\end{split}
\end{equation}
We intend to show that with $k\to\infty$ this implies \eqref{QPP5} and hence, completes the proof.

By It\^o's formula and using \eqref{SPP1} and \eqref{SPPn}, we may find out  a constant $C>0$ such that
\beq\label{MM3} \beg{split} \d |X_t-X_t^{(k)}|^2\le \d M_t + &C\big\{|X_t-X_t^{(k)}|^2 + \E |X_t-X_t^{(k)}|^2+\\
& \qquad \|\si(t,X_t,\mu_t)- \si^{(k)}(t,X_t,\mu_t)\|_{HS}^2\big\}\d t,~~~ t\in [0,T] \end{split}
\end{equation}
holds for some martingale $M_t$. By \eqref{SPP1} and the definition of $\si^{(k)}$, for some constant $C'>0$ we have
$$\|\si(t,X_t,\mu_t)- \si^{(k)}(t,X_t,\mu_t)\|_{HS}^2\le C' |X_t-\phi_k(X_t)|^2\le C'|X_t|^21_{\{|X_t|\ge k\}}.$$
Combining this with \eqref{MM3} and using Gronwall's lemma, we arrive at
$$\E |X_t-X_t^{(k)}|^2\le CC'\e^{2Ct} \int_0^t\E\big[|X_s|^21_{|X_s|\ge k}\big]ds,\ \ t\in [0,T].$$
This, together with \eqref{MM}, implies
$$\lim_{k\to\infty} \sup_{t\in [0,T]} \\W_2(\mu_t,\mu_t^{(k)})^2 \le \lim_{k\to\infty} \sup_{t\in [0,T]} \E |X_t-X_t^{(k)}|^2=0.$$
In particular, $\{\mu_t^{(k)}:\ t\in [0,T], k\ge 1\}$ is compact in $\scr P_2(\R^d).$ So,
from  the continuity of $\si b, \pp_\mu f$, and $\pp_y\pp_\mu f,$ the linear growth of $|\si b|$, and the condition   $f\in \scr C([0,\infty)\times\R^d,\scr P_2(\R^d))$,
it is easy to see that   with $k\to\infty$   \eqref{QPP5'} implies  \eqref{QPP5}.
\end{proof}

We will also need the following result which is more or less standard for classical SDEs. For readers' convenience we present a complete proof for the present distribution dependent setting.

\beg{lem}\label{LNN} Assume $(H)$. For any $s\ge 0$ and $\mu\in \scr L_2(\R^d)$, let $(X_t)_{t\ge s}$ solve $\eqref{E1}$ for $\L_{X_s}=\mu$. Then
for any $t>s$, $\L_{X_t}$ has has full support on $\R^d$. \end{lem}
\beg{proof}  Consider the SDE \eqref{E1'*}, we have
\beq\label{LNN2} \L_{X_t}= \int_{\R^d} \L_{X_{s,t}^{x,\mu}}\mu(\d x),\ \ t>s.\end{equation}So, it suffices to prove that for any $x\in\R^d$ and $t>s$, $\L_{X_{s,t}^{x,\mu}}$ has full support on $\R^d$. Let
$$P_{s,t}^\mu f(x)= \E f(X_{s,t}^{x,\mu}),\ \ t>s, f\in \B_b(\R^d).$$
By \cite[Theorem 5.1]{LAP}, assumption $(H)$ implies the   Harnack inequality  
\beq\label{LHH} P_{s,r}^\mu f(y)\le \psi(r,t,x,y)P_{s,t}^\mu f(x),\ \ t>r>s, x,y\in \R^d, 0\le f\in \B_b(\R^d) \end{equation}
for some  measurable function $$\psi: (r,t,x,y)\mapsto (1,\infty),\ \ t>r>s, x,y\in \R^d.$$ 
If for some $t>s$ and $x\in\R^d$ the distribution $\L_{X_{s,t}^{x,\mu}}$ does not have full support $\R^d$, then there exist $y\in\R^d$ and $\vv>0$ such that $$P_{s,t}^\mu 1_{B(y,\vv)}(x)=\P(X_{s,t}^{x,\mu}\notin B(y,\vv))=0.$$
Combining this with \eqref{LHH} gives
$$\P(X_{s,r}^{y,\mu}\in B(y,\vv))= P_{s,r}^{\mu} 1_{B(y,\vv)}(y) \le \psi(r,t,x,y)P_{s,t}^\mu 1_{B(y,\vv)}(x) =0,\ \ r\in (s,t).$$
By the continuity of $X_{s,r}^{y,\mu}$, by letting $r\to s$ we obtain $\P(y=X_{s,s}^{y,\mu}\in B(y,\vv/2))=0$ which is impossible. 
So, as required, for any $x\in\R^d$ and $t>s$ the distribution $\L_{X_{s,t}^{x,\mu}}$ has full support on $\R^d$.\end{proof}

\beg{proof}[Proof of Theorem \ref{T1.1}] (1)
Let $\mu_t= \L_{X_t}$. Applying the It\^o formula \eqref{rpp}  yields
\begin{equation}\label{YUU}
\begin{split}
\d V(t,X_t, \mu_t)&= (\pp_t V+ {\bf L}_{\si,b} V)(t,X_t, \mu_t) \d t  +\<(\si^*\pp_x V)(t,X_t, \mu_t), \d W_t\>.\\
\end{split}
\end{equation}
This, together with \eqref{ATT0}, gives
\beq\label{YUU1}
\d V(t, X_t, \mu_t)={\bf f}(t,X_t,\mu_t)\d
t+\<{\bf g}(t,X_t,\mu_t),\d W_t\>,~~~~t\geqslant 0.
\end{equation}
Whence, \eqref{YPP0} follows by integrating \eqref{YUU1} from $s$ to $t.$

(2) On the other hand, for any $s\in [0,T)$ and $\mu\in \scr P_2(\R^d)$, let $X_s\in L^2(\OO\to\R^d,\F_t,\P)$ with $\scr L_{X_s}=\mu$, and let $(X_t)_{t\in [s, T]}$ solve \eqref{E1} from time $s$.  By combining \eqref{YUU} with \eqref{YUU1} and using the uniqueness of decomposition for semi-martingale, we infer that
\begin{equation*}
{\bf f}(t,X_t, \mu_t)=(\pp_t V+ {\bf L}_{\si,b} V)(t,X_t, \mu_t) ,~~~~{\bf g}(t,X_t, \mu_t)=(\si^*\pp_x V)(t,X_t, \mu_t),\ \ t\in [s,T],
\end{equation*} where $\mu_t:=\L_{X_t}$ with $\mu_s=\mu$.
Since by Lemma \ref{LNN} the assumption $(H)$ implies that $\mu_t$ for $t\in (s,T]$  has a full support on $\R^d$,  we derive
\begin{equation*}
{\bf f}(t,x, \mu_t)=(\pp_t V+ {\bf L}_{\si,b} V)(t,x, \mu_t) ,~~~~{\bf g}(t,x, \mu_t)=(\si^*\pp_x V)(t,x, \mu_t),\ \   x\in \R^d,  t\in (s,T].
\end{equation*} Since $\mu_t$ is continuous in $t$,  and since $f, (\pp_t +{\bf L}_{\si,b})V $ are continuous on $[0,T]\times\R^d\times \scr P_2(\R^2)$, by letting $t\downarrow s$ we obtain
$${\bf f}(s, x, \mu)=(\pp_s V+ {\bf L}_{\si,b} V)(s,x, \mu) ,~~~~{\bf g}(s,x, \mu)=(\si^*\pp_x V)(s,x, \mu),\ \  x\in \R^d. $$
By the arbitrariness  of $s\in [0,T)$ and $\mu\in \scr P_2(\R^d)$, this implies  \eqref{ATT0}.
\end{proof}

To prove Theorem \ref{T1.2}, we will need the following lemma, which reduces to  the main result Theorem 6.2  in \cite{LP} when $b(t,x,\mu)$ and $\si(t,x,\mu)$ are independent of $t$. Since the proof of \cite[Theorem 6.2]{LP} also applies  to the the present time inhomogeneous situation,  we skip the proof.

\beg{lem}[\cite{LP}]\label{LLP} In the situation of Theorem $\ref{T1.2}$, let  $\Phi\in C_b^{2,(1,1)}(\R^d\times\scr P_2(\R^d))$.
Then $V(t,x,\mu):= \E\Phi(X_{t,T}^{x,\mu}, P_{t,T}^*\mu) $ is the unique solution to the PDE
\begin{equation}\label{YGG}
\begin{cases}
(\pp_t+\mathbf L_{\si,b})V(t,x,\mu)=0,\\
 V(T,x,\mu)=\Phi(x,\mu),\ \  \ ~~~ (t,x,\mu)\in [0,T]\times\R^d\times\scr P_2(\R^d).
\end{cases}
\end{equation}\end{lem}

We will also need the following lemma for a probabilistic representation of a particular solution  to  the first equation in \eqref{ATT0}.

\beg{lem}\label{LLP2} In the situation of Theorem $\ref{T1.2}$, let
$$V_\f(t,x,\mu)=  -\E\int_t^T\f(r,X_{t,r}^{x,\mu}, P_{t,r}^*\mu)\d r,\ \ (t,x,\mu)\in [0,T]\times \R^d\times\scr P_2(\R^d).$$
Then $V_\f$ is the unique solution to the PDE
\begin{equation}\label{YGG1}
\begin{cases}
(\pp_t+\mathbf L_{\si,b} ) V_{\f}(t,x,\mu)=\f(t,x,\mu), \\  V_{\f}(T,x,\mu)=0,\ \   (t,x,\mu)\in [0,T]\times \R^d\times\scr P_2(\R^d).
\end{cases}
\end{equation}
 \end{lem}

\beg{proof}
(a) We first observe that $ V_{\f}(t,x,\mu)$ solves \eqref{YGG1}. Obviously,
$$ V_{\f}(T,x,\mu)=0.$$
It remains to prove
\begin{equation}\label{GGMM1}
(\pp_t+\mathbf L_{\si,b} ) V_{\f}(t,x,\mu)=\f(t,x,\mu).
\end{equation}
By the definition of $ V_{\f}$ and our condition on $\f$, we have
$$(\mathbf L_{\si,b} V_{\f})(t,x,\mu)=-\int_t^T \mathbf L_{\si,b}\{\E \f(r,X_{t,r}^{x,\mu}, P_{t,r}^*\mu)\}\d r,$$
and
$$(\pp_t V_{\f})(t,x,\mu)=\f(t,x,\mu)-\int_t^T \pp_t\{\E \f(r,X_{t,r}^{x,\mu}, P_{t,r}^*\mu)\}\d r.$$
So,
$$(\pp_t+\mathbf L_{\si,b})V_{\f}(t,x,\mu)=\f(t,x,\mu)- \int_t^T(\pp_t+ \mathbf L_{\si,b})\{\E \f(r,X_{t,r}^{x,\mu}, P_{t,r}^*\mu)\}\d r.$$
On the other hand, applying Lemma \ref{LLP} to $T=r$ and $\Phi(x,\mu)=\f(r,x,\mu)$, we obtain
$$(\pp_t+ \mathbf L_{\si,b})\{\E \f(r,X_{t,r}^{x,\mu}, P_{t,r}^*\mu)\}=0,~~~~~~r\in(t,T]. $$
Therefore, \eqref{GGMM1} holds.

(b) We assume that $U(t,x,\mu)$ is another solution to \eqref{YGG1} with $U(T,x,\mu)=0$. By Lemma \ref{L1},  for any $0\leqslant t \leqslant s \leqslant T$,
$$U(s,X_{t,s}^{x,\mu}, P_{t,s}^*\mu)-\int_{t}^{s} \f(u,X_{t,u}^{x,\mu}, P_{t,u}^*\mu)\d u$$
and $$V_{\f}(s,X_{t,s}^{x,\mu}, P_{t,s}^*\mu)-\int_{t}^{s} \f(u,X_{t,u}^{x,\mu}, P_{t,u}^*\mu)\d u$$
are martingales. Then
$$U(s,X_{t,s}^{x,\mu}, P_{t,s}^*\mu)-V_{\f}(s,X_{t,s}^{x,\mu}, P_{t,s}^*\mu)$$
is a martingale. Combining this with  $U(T,x,\mu)=V_{\f}(T,x,\mu)=0$, we arrive at
$$U(t,x,\mu)-V_{\f}(t,x,\mu)=\E \big(U(T,X_{t,T}^{x,\mu}, P_{t,T}^*\mu)-V_{\f}(T,X_{t,T}^{x,\mu}, P_{t,T}^*\mu)|\F_t\big)=0.$$
Then the uniqueness is proved.
\end{proof}

\beg{proof}[Proof of Theorem \ref{T1.2}] By Theorem \ref{T1.1}, it suffices to prove the first assertion.
By Lemma \ref{LLP}, we deduce that
$$V_1(t,x,\mu):= \E\Phi(X_{t,T}^{x,\mu}, P_{t,T}^*\mu)$$
is the unique solution to the PDE \eqref{YGG}.
And, according to Lemma \ref{LLP2}, we know that
$$V_{\f}(t,x,\mu):=  -\E\int_t^T\f(r,X_{t,r}^{x,\mu}, P_{t,r}^*\mu)\d r$$ solves \eqref{YPP1}.
So,
\begin{equation*}
V(t,x,\mu):=V_1(t,x,\mu)+V_{\f}(t,x,\mu)=\E\Phi(X_{t,T}^{x,\mu}, P_{t,T}^*\mu)-\E\int_t^T\f(r,X_{t,r}^{x,\mu}, P_{t,r}^*\mu)\d r
\end{equation*}
together with $\g=(\si^*\pp_tV)(t,x,\mu)$, solves \eqref{ATT0}.

On the other hand, let $V(t,x,\mu)$ solve \eqref{ATT0} and let
$$\Phi(x,\mu)= V(T,x,\mu),\ \ (x,\mu)\in \R^d\times\scr P_2(\R^d).$$ It suffices to prove
\begin{equation}\label{GGMM}
V(t,x,\mu)=\E \bigg(\Phi(X_{t,T}^{x,\mu}, P_{t,T}^*\mu)-\E\int_t^T\f(r,X_{t,r}^{x,\mu}, P_{t,r}^*\mu)\d r\bigg).\end{equation}
 Indeed, by \eqref{ATT0} and Lemma \ref{LLP2}, we have
$$(\pp_t+\mathbf L_{\si,b} )( V-V_{\f})(t,x,\mu)=0.$$So, Lemma \ref{LLP} and $V_{\f}(T,x,\mu)=0$ imply
$$( V-V_{\f})(t,x,\mu)=\E\Phi(X_{t,T}^{x,\mu}, P_{t,T}^*\mu) $$
with  $(V-V_{\f})(T,x,\mu)=V(T,x,\mu) =\Phi(x,\mu)$. This, together with the definition of $V_{\f}$, implies \eqref{GGMM}. Then the proof is completed.
\end{proof}

\beg{proof}[Proof of Corollary \ref{T1.4}]  Assertion (1) is direct consequence of Theorem \ref{T1.1} for $\f= \ff 1 {2\bb} |\g|^2$. It remains to prove assertion (2).

Under the condition of assertion (2), let $\tt V(t,x,\mu)= \E \Phi(X_{t,T}^{x,\mu}, P_{t,T}^*\mu).$ By Lemma \ref{LLP}  we have
$$(\pp_t +{\bf L}_{\si,b})\tt V(t,x,\mu)=0.$$
Since for $V$ in \eqref{TTY0}  we have $V= -\bb\log\tt V$, this implies
\beg{align*}  (\pp_t +{\bf L}_{\si,b}) V(t,x,\mu)
& = -  \ff {\bb (\pp_t + {\bf L}_{\si,b}) \tt V } {\tt V}(t,x,\mu) + \ff{\bb|\si^*\pp_x\tt V|^2(t,x,\mu)}{2 \tt V^2(t,x,\mu)} \\
 &=  \ff 1 {2\bb} |\si^*\pp_x  V|^2.\end{align*} So, \eqref{NPDE} holds, and the boundary condition   $V(T,x,\mu)= -\bb \log\Phi(x,\mu)$ follows from \eqref{TTY0} and the definition  of  $\tt V$.

On the other hand, let $V\in C_b^{1,2,(1,1)}([0,T]\times\R^d\times\scr P_2(\R^d))$ solve \eqref{NPDE}. We take
\begin{equation}\label{GMA}
\tt V(t,x,\mu)=\exp[-\beta ^{-1}V(t,x,\mu)],\ \ (t,x,\mu)\in [0,T]\times\R^d\times\scr P_2(\R^d).
\end{equation}
It is easy to see that \eqref{NPDE} implies
$$(\pp_t+ {\bf L}_{\si,b})\tt V (t,x,\mu)=0,~~~ (t,x,\mu)\in [0,T]\times \R^d\times \scr P_2(\R^d).$$
Therefore, by Lemma \ref{LLP}   we have
$$
\tt V(t,x,\mu)=\E \tt V (T,\tt X_{t,T}^{x,\mu}, \tt\mu_{t,T})=:\E\Phi(\tt X_{t,T}^{x,\mu}, \tt\mu_{t,T})~~~(t,x,\mu)\in [0,T]\times \R^d\times \scr P_2(\R^d).
$$
Combining this with \eqref{GMA}, we obtain \eqref{TTY0} and hence finish the proof.
\end{proof}

\beg{proof}[Proof of Corollary \ref{C1.3}]
By \eqref{XTT} and \eqref{SIG}, the definitions of $\mathbf L_{\si}$ and $\mathbf L_{\si,b}$, we have
 \begin{equation*}
 \begin{split}
 &(\pp_t+ {\bf L}_{\si}) V (t,x,\mu)\\
 &= \pp_tV (t,x,\mu)+\ff 1 2 {\rm tr}\big(\si\si^*  \pp_{x}^2 V)(t,x,\mu) +\big\<b,\pp_xV\big\> (t,y,\mu) \\
 &\quad+\int_{\R^d} \Big[\ff 1 2 {\rm tr}\big\{(\si\si^*)(t,y,\mu)  \pp_{y}\pp_\mu V (t,x,\mu)(y)\big\} +\big\< b(t,y,\mu) ,  \pp_{\mu} V(t,x,\mu)(y)\big\>\Big]\mu(\d y)\\
 &\quad+\ff 1 2|\si^* \pp_x V|^2  (t,x,\mu)-\big\<b,\pp_xV\big\> (t,y,\mu) \\
 &\quad+\int_{\R^d} \big\<(\si\si^*\pp_yV )(t,y,\mu) ,  \pp_{\mu} V(t,x,\mu)(y)\big\>\mu(\d y)\\
 &\quad -\int_{\R^d} \big\<b(t,y,\mu) ,  \pp_{\mu} V(t,x,\mu)(y)\big\>\mu(\d y)\\
 &=(\pp_t+ {\bf L}_{\si,b}) V (t,x,\mu)+\ff 1 2|\si^* \pp_x V|^2  (t,x,\mu)-\big\<b,\pp_xV\big\> (t,y,\mu) \\
 &\quad +\int_{\R^d} \big\<(\si\si^*\pp_yV-b)(t,y,\mu) ,  \pp_{\mu} V(t,x,\mu)(y)\big\>\mu(\d y).
 \end{split}
 \end{equation*}
 Combining this with $b(t,x,\mu)= (\si\si^*\pp_x V)(t,x,\mu)$, we obtain
 \begin{equation}\label{NPY}
 (\pp_t+ {\bf L}_{\si}) V (t,x,\mu)=(\pp_t+ {\bf L}_{\si,b}) V (t,x,\mu)-\frac{1}{2}|\si^*\pp_xV|^2(t,x,\mu).
\end{equation}
We are now ready to finish the proof by using Theorem \ref{T1.1} and Corollary \ref{T1.4}.

If \eqref{ATT2} holds, then  \eqref{BB} holds for $\tt b=\si^*\pp_x V$, and  \eqref{NPY} implies \eqref{ATT0} for $\f(t,x,\mu)=\frac{1}{2}|\tt b|^2(t,x,\mu)$ and $\g(t,x,\mu)=\tt b(t,x,\mu)$. So, by Theorem \ref{T1.1}(1), $A_{s,t}$ is path independent.
On the other hand, if  \eqref{BB} holds for $\tt b=\si^*\pp_x V$ and $A_{s,t}$ is path independent in the sense of \eqref{YPP0} for some $V\in C^{1,2,(1,1)}([0,T]\times\R^d\times\scr P_2(\R^d))$, then  by Theorems \ref{T1.1}(2) and \eqref{NPY},   \eqref{ATT2} holds. So, assertion (1) is proved.

Finally, by \eqref{NPY}, the first equation in \eqref{ATT2} is equivalent to \eqref{NPDE} for $\bb=1$. Then the second assertion (2) follows from Corollary \ref{T1.4}(2) for $\bb=1.$
 \end{proof}


\end{document}